\documentclass[a4paper,11pt]{article}
\usepackage{amsfonts,amsmath,amssymb}
\usepackage[a4paper]{geometry}
\usepackage[pdftex]{graphicx}

\newcommand{\vect}[1]{\mathbf{#1}}

\newcommand{\R}{\mathbb{R}}
\newcommand{\kron}{\otimes}
\newcommand{\Real}{\Re\mathfrak{e}}
\newcommand{\Imag}{\Im\mathfrak{m}}
\newcommand{\factorm}{q_j\left(s,\,P,\,D,\,\tau_1,\,\tau_2,\,\lambda_j\right)}

\newtheorem{lemma}{Lemma}
\newtheorem{definition}{Definition}
\newtheorem{proposition}{Proposition}

\begin{document}
\title{A Consensus Protocol under Directed Communications with Two Time Delays and Delay Scheduling}

\author{Rudy Cepeda-Gomez\thanks{Faculty of Mechatronic Engineering, Universidad Santo Tomas, Bucaramanga, Colombia. Email: \mbox{rudycepedagomez@mail.ustabuca.edu.co}} and Nejat Olgac\thanks{Mechanical Engineering Department, University of Connecticut, Storrs, CT, USA. Email: \mbox{olgac@engr.uconn.edu}}}
\date{}

\maketitle

\begin{abstract}
This paper studies a consensus protocol over a group of agents driven by second order dynamics. The communication among members of the group is assumed to be directed and affected by two rationally independent time delays, one in the position and the other in the velocity information channels. These delays are unknown but considered to be constant and uniform throughout the system. The stability of the consensus protocol is studied using a simplifying factorization procedure and deploying the Cluster Treatment of Characteristic Roots (CTCR) paradigm.  This effort results in a unique depiction of the exact stability boundaries in the domain of the delays. CTCR requires the knowledge of the potential stability switching loci exhaustively within this domain. The creation of these loci is an important contribution of this work. It is done in a new surrogate coordinate system, called the \emph{Spectral Delay Space (SDS)}. The relative stability, i.e., the speed to reach consensus of the system is also investigated for this class of systems. Based on the outcome of this effort, a paradoxical control design concept is introduced. It is called the \emph{Delay Scheduling}, which is another key contribution of this paper. It reveals that the performance of the system may be improved by increasing the delays. The amount of increase, however, is only revealed by the CTCR. Example case studies are presented to verify the underlying analytical derivations.
\end{abstract}

\section{Introduction}
Recently considerable research effort has been placed in the multi-agent systems area, and particularly in the consensus behaviour of swarms. The consensus aims at an agreement among the agents. This problem is treated by \cite{vicsek1995} in the early on, where agents aim to align their headings using discrete time representation. Later, \cite{olfati-saber2004} introduced a formal framework of the consensus problem for agents driven by first order dynamics. That study considers directed communication topologies, both for fixed and switching types, and introduce some graph theoretical results useful for the stability analysis of such agreement protocols. Under the simplifying features of first order dynamics, they also consider time delayed communications in the case of a fixed topology. Many other researchers suggest further extensions to this work, proposing consensus protocols for agents driven by second order dynamics \cite{ren2005,ren2007} without time delays in the communication channels. 

The works by Lin and his co-workers include the study of agents driven by second order discrete \cite{Lin2009-2} and continuous \cite{Lin2009-1} dynamics, and applications of different control laws, ranging from a basic $P$-$D$-like (proportional and integral) control logic \cite{Lin2008-1}, to $H_\infty$ structure \cite{Lin2008-2}. Most of these works entail time delayed communications and directed topologies. Sun and Wang also consider the consensus problem with fixed \cite{Sun2009-1} and time varying delays \cite{Sun2009-2}. \cite{Peng2009} study a leader follower case, also with time varying delays. \cite{Meng2011} presents a comparison study among several different swarm control laws is, including first and second order agents and considering two time delays. 

When it comes to the stability analysis of consensus systems with time delays, almost all the previous works rely on Lyapunov-Krasovskii or Razhumikin based methodologies \cite{Lin2009-2,Lin2009-1,Lin2008-1,Lin2008-2,Meng2011,olfati-saber2004}. They produce only very conservative results, leading to stability bounds for very small delays. Additionally, because these results are based on some solutions of an LMI, they are always imprecise, conservative and, very critically, cumbersome to deploy.

As the non-delayed control structure imparts a stable consensus, conservatively small delays would still maintain the stable posture. Very differently from the general trend, we wish to determine the stability bounds in delay non-conservatively and exhaustively in this paper. This venue opens up an unusual path of control logic design, which is called the \emph{Delay Scheduling}.

Recently the authors have presented a methodology for the analysis of consensus protocols with single and multiple, rationally-independent time delays \cite{TAC}. This earlier work is based on the combination of a simplifying factorization procedure over the characteristic equation of the system and the deployment of a crucial stability paradigm, which is called the Cluster Treatment of Characteristic Roots, CTCR, \cite{Olgac2005,Fazelinia2007,Ergenc2007}. CTCR provides a tool for assessing the stability of LTI systems with multiple rationally-independent delays. This new method uniquely creates exact, exhaustive and explicit stability regions in the domain of the delays. All of our earlier investigations \cite{TAC} have focused on undirected communications topologies. In the present document, on the other hand, we consider directed topologies. The transition from undirected to directed topologies is not trivial and the treatment of the ensuing difference forms the critical contribution of this document. Although some other investigators have treated this problem, as mentioned before, this is the first work, to the best of our knowledge, in which truly large time delays and exact stability boundaries are handled precisely (non-conservatively) and exhaustively. The "large" delay, here, implies that they are either orders of magnitude larger than the control sampling period or they are in the order of the fastest desired dynamics in the system. 

In summary, the two main contributions of this paper are as follows:  for a group dynamics to achieve consensus under time-delayed communication topologies:
\begin{enumerate}
\item	The \emph{delay scheduling} concept is introduced for swarm control for the first time. It is a direct consequence of the stability assessment paradigm, called CTCR, which declares the exact and exhaustive stable regions in the domain of the two rationally independent time delays.
\item	This treatment is deployed for the first time for directed inter-agent communication topologies on a group dynamics with second order governing equations.
\end{enumerate}

Bold face notation is used in the text for vector quantities, bold capital letters for matrices and italic symbols for scalars.
\section{Consensus Protocol}
We consider a group of $n$ autonomous agents, which are driven by second order dynamics given by $\ddot{x}_j\left(t\right)=u_j\left(t\right)$, $j=1,\,2,\ldots,\,n$, where $x_j\left(t\right)\in\R$ is taken as the scalar position and $u_j\left(t\right)\in\R$ as the control law, without loss of generality. The arguments are presented for a one-dimensional case for simplicity, and all the results are scalable to higher dimensions by using the Kronecker product notation, the additional dimensions are decoupled in the governing dynamics. Features such as "holonomicity" and "non-linearity" in the agent dynamics would also prevent this expansion.

Following the convention in the literature \cite{olfati-saber2004,Lin2008-1,Lin2008-2,Lin2009-1,Meng2011,Sun2009-1,Sun2009-2,Peng2009,TAC}, we declare the consensus is achieved when all $n$ agents are at the same position, i.e., $\lim_{t\rightarrow\infty}\left(x_j\left(t\right)-x_k\left(t\right)\right)=0$ for any $j,k\in\left[1,\,n\right]$. Notice that this consensus definition does not state anything \emph{a priori} regarding the value of the final position.
 
The $j$-th agent receives position and velocity information from a subset of agents, $\mathcal{N}_j$, which consists of $\delta_j$ agents, $\delta_j<n$, and it is called the \emph{informers} of the $j$-th agent. Assuming unidirectional channels, the communication network can be described by a directed graph with $n$ vertices, and $\delta_j$ is the in-degree of node $j$. It is also assumed that all these communication channels have two constant and uniform time delays, $\tau_1$ and $\tau_2$, which affect the position and velocity information exchange, respectively.  
In order to create the consensus we propose a $PD$ (proportional and derivative) decentralized control law for each agent $j$:
\begin{equation}
\ddot{x}_j\left(t\right)=P\left(\sum_{k\in\mathcal{N}_j}{\frac{x_k\left(t-\tau_1\right)}{\delta_j}}-x_j\left(t\right)\right)+D\left(\sum_{k\in\mathcal{N}_j}{\frac{\dot{x}_k\left(t-\tau_2\right)}{\delta_j}}-\dot{x}_j\left(t\right)\right)
\label{eq:c2delay1}
\end{equation}
where $P$ and $D$ are user-selected, positive control gains. This logic tries to bring the agent's current position to the centroid of its informers, $\mathcal{N}_j$, and its velocity to the mean velocity of the same set of informers (or the velocity of the same centroid), using the last known positions and velocities of the informer agents. Notice that this protocol does not include self delayed information, while all the data from the informers is delayed, positions by $\tau_1$ and velocities by $\tau_2$. The corresponding dynamics of the complete system in state space becomes:
\begin{equation}
\begin{split}
\dot{\vect{x}}\left(t\right)=&\left(\vect{I}_n\kron\left[\begin{array}{cc}0&1\\-P&-D\end{array}\right]\right)\vect{x}\left(t\right)\\&+\left(\vect{C}\kron\left[\begin{array}{ccc}0&0\\P&0\end{array}\right]\right)\vect{x}\left(t-\tau_1\right)+\left(\vect{C}\kron\left[\begin{array}{ccc}0&0\\0&D\end{array}\right]\right)\vect{x}\left(t-\tau_2\right)
\end{split}
\label{eq:c2delay2}
\end{equation}
with the state vector being  a concatenation of the $n$ positions and velocities: $\vect{x}=\left[x_1\ \dot{x}_1\ x_2\ \dot{x}_2\  \cdots\ x_n\ \dot{x}_n\right]^T\in\R^{2n}$. In \eqref{eq:c2delay2}, $\kron$ denotes Kronecker multiplication \cite[see]{kron}, $\vect{I}_n$ is the identity matrix of dimension $n$ and $\vect{C}=\boldsymbol{\Delta}^{-1}\vect{A}_\Gamma$ is the weighted adjacency matrix, created by the multiplication of the inverse of $\boldsymbol{\Delta}=\text{diag}\left(\delta_1,\,\delta_2,\ldots,\,\delta_n\right)\in\R^{n\times n}$, known as the in-degree matrix, and $\vect{A}_{\Gamma}$, the adjacency matrix of the communication topology, following graph theory conventions \cite[see]{Godsil,Biggs}. 

The dynamics in \eqref{eq:c2delay2} can be expressed in a more compact form as:
\begin{equation}
\dot{\vect{x}}\left(t\right)=\vect{Ax}\left(t\right)+\vect{B}_1\vect{x}\left(t-\tau_1\right)+\vect{B}_2\vect{x}\left(t-\tau_2\right)
\label{eq:genmtds}
\end{equation}
with self evident $\vect{A}$, $\vect{B}_1$ and $\vect{B}_2$ matrices from \eqref{eq:c2delay2}.

The complexity level of this dynamics increases rapidly as the number of agents gets larger, making the stability analysis numerically tedious and intractable. The corresponding characteristic equation of the system is
\begin{equation}
Q\left(s,\,P,\,D,\,\tau_1,\,\tau_2\right)=\det\left(s\vect{I}_{2n}-\vect{A}-\vect{B}_2e^{-\tau_1\,s}-\vect{B}_2e^{-\tau_2\,s}\right)=0
\label{eq:chareq}
\end{equation}
This equation is a $2n$ order quasi-polynomial with commensuracy degree $n$ (i.e., up to $n\tau_1$ and $n\tau_2$ terms appear) and, more critically, the cross-talk between the time delays (i.e., terms like $\tau_1+\tau_2$) are typically present. The complexity of \eqref{eq:chareq} increases further if a more tortuous communication topology is used (e.g., one with non-uniform delay structure). Such general stability problem of multiple time delay systems is notoriously known to be NP-hard \cite{Olgac2005,Toker1996}. Even the strongest of the existing mathematical tools to analyse the stability of time delay systems of this class falls short to handle such cases exhaustively. Two important procedures are combined to remedy the impasse $vis-\grave{a}-vis$ the particular consensus problem with directed topologies. First, a factorization operation, introduced by \cite{TAC}, is performed to break the characteristic equation into factors of reduced order quasi-polynomial and simpler but identical forms. The crucial difference of this paper from our earlier work appears here. The directed topologies result in complex conjugate identifiers of the factors, as opposed to all real in the undirected cases. Second, a method, called the Cluster Treatment of Characteristic Roots (CTCR) \cite{Olgac2005,Fazelinia2007,Ergenc2007} and a new domain, Spectral Delay Space (SDS), which we describe in Section \ref{sec:III}, are deployed for the first time on directed topologies to perform and exact and exhaustive stability analysis. This combination of handling complex eigenvalues and SDS facilitates an efficient and novel mechanism to resolve the stability problem at hand which constitutes the main contribution of the present work.
\subsection{Factorization Procedure}
The stem procedure in \cite{TAC}, which transforms the characteristic equation \eqref{eq:chareq} into a product of a set of second order factors, was based on the use of an undirected communication structure. For undirected case, the eigenvalues of the matrix $\vect{C}$ are always real. When directed topologies are considered, however, the matrix $\vect{C}$ may also have complex conjugate eigenvalues. While this does not affect the main philosophy behind the methodology we present, the practical implementation changes substantially.

Assume that there is a non-singular matrix $\vect{T}$ which converts $\vect{C}$ into its Jordan canonical form: $\boldsymbol{\Lambda}=\vect{T}^{-1}\vect{CT}$. The matrix $\boldsymbol{\Lambda}$ is block diagonal of the form:
\begin{equation}
\boldsymbol{\Lambda}=\left[\begin{array}{ccccccc}
\lambda_1&0&\cdots&0&0&\cdots&0\\
0&\lambda_2&\cdots&0&0&\cdots&0\\
\vdots&\vdots&\ddots&\vdots&\vdots&\ddots&\vdots\\
0&0&\cdots&\lambda_{\ell}&0&\cdots&0\\
\vect{0}_{21}&\vect{0}_{21}&\cdots&\vect{0}_{21}&\vect{J}_{\ell+1}&\cdots&\vect{0}_{21}\\
\vdots&\vdots&\ddots&\vdots&\vdots&\ddots&\vdots\\
\vect{0}_{21}&\vect{0}_{21}&\cdots&\vect{0}_{21}&\vect{0}_{21}&\cdots&\vect{J}_{\ell+m}\end{array}\right]
\label{eq:lambdadirect}
\end{equation}
where $\vect{0}_{21}=\left[0\ 0\right]^T$; $\lambda_j$, $j=1,\,2\,\ldots,\,n$, are the (size 1) Jordan blocks corresponding to the real eigenvalues of $\vect{C}$ and 
\begin{equation}
\vect{J}_j=\left[\begin{array}{cc}
\Real\left(\lambda_j\right)&-\Imag\left(\lambda_j\right)\\
\Imag\left(\lambda_j\right)&\Real\left(\lambda_j\right)\end{array}\right],\qquad j=\ell+1,\,\ell+2,\,\ldots,\,\ell+m
\label{eq:jordanblock}
\end{equation}
are the $2\times2$ Jordan blocks corresponding to the complex conjugate eigenvalue pairs.

\emph{Remark I}. It is trivial to observe that transformation $\vect{x}\left(t\right)=\left(\vect{I}_2\kron\vect{T}^{-1}\right)\boldsymbol{\xi}$ converts \eqref{eq:c2delay2} into $\ell+m$ decoupled subsystems, which generate the factors of the characteristic equation, as we display next. The subsystems corresponding to real eigenvalues (size 1 Jordan blocks) create factors of the form:
\begin{subequations}
\label{eq:factorsdir}
\begin{equation}
\factorm=s^2+Ds+P-\lambda_j\left(Dse^{-\tau_2s}+Pe^{-\tau_1s}\right)
\label{eq:2ndfactor}
\end{equation}
whereas the size two Jordan block corresponding to complex conjugate eigenvalue pairs generate factors like:
\begin{equation}
\begin{split}
q_j&\left(s,\,P,\,D,\,\tau_1,\,\tau_2,\,\lambda_j\right)=s^4+2Ds^3+\left(D^2+2P\right)s^2+2DPs+P^2\\
&-2\Real\left(\lambda_j\right)\left(s^2+Ds+P\right)\left(Dse^{-\tau_2s}+Pe^{-\tau_1s}\right)+\left|\lambda_j\right|^2\left(Dse^{-\tau_2s}+Pe^{-\tau_1s}\right)^2
\end{split}
\label{eq:4thfactor}
\end{equation}
\end{subequations}

The characteristic equation of the complete system is then formed by the product of these $\ell+m$ individual factors.  
\begin{equation}
\begin{split}
&Q\left(s,\,P,\,D,\,\tau_1,\tau_2\right)=\prod_{j=1}^{\ell+m}\factorm=\\
&\prod_{j=1}^{\ell}\left[s^2+Ds+P-\lambda_j\left(Dse^{-\tau_2s}+Pe^{-\tau_2s}\right)\right]=\prod_{j=\ell+1}^{\ell+m}\left[s^4+2Ds^3+\left(D^2+2P\right)s^2\right.\\
&\quad\left.+2DPs+P^2-2\Real\left(\lambda_j\right)\left(s^2+Ds+P\right)\left(Dse^{-\tau_2s}+Pe^{-\tau_2s}\right)\right.\\
&\quad\left.+\left|\lambda_j\right|^2\left(Dse^{-\tau_2s}+Pe^{-\tau_2s}\right)^2\right]=0
\end{split}
\label{eq:cefactdir}
\end{equation}
Clearly, this procedure simplifies the problem considerably, by transforming it from a $2n$ order system with time delays of commensuracy degree $n$ and delay cross-talk  into a combination of $\ell$ second order and $m$ fourth order subsystems with maximum commensuracy degree of 2 (e.g., $e^{-2\tau_1\,s}$) and single delay cross-talk (e.g., $e^{-\left(\tau_1+\tau_2\right)s}$). Furthermore, since the only differentiating element from one factor to the other is the eigenvalue $\lambda_j$, the stability analysis in the domain of the delays can be performed just once in each case, for a generic real $\lambda$ and a generic complex $\lambda$. The successive repetitions of the results from this for $\ell+m$ different subsystems will suffice to obtain the stability picture of the complete system. This feature provides us a very strong computational platform for the stability analysis.

\emph{Remark II}. We now wish to direct the discussions to some other features of the eigenvalues of the weighted adjacency matrix $\vect{C}$. From the way in which this matrix is created, its diagonal consists of zeros only, and any row elements always add up to 1. The latter property makes $\vect{C}$ a row-stochastic matrix of which the components are all non-negative \cite{Marcus1996}. Row-stochastic non-negative matrices possess a wonderful feature: the eigenvalues of such matrices are always within the unit disk of the complex plane. This feature arises as an extension to the Gershgorin's disk theorem \cite{Bell1965}. 

\emph{Remark III}. Furthermore, it is proven by \cite{Agaev2005} that if the topology is connected and has at least one spanning tree, $\lambda=1$ is one of the eigenvalues of the weighted adjacency matrix, $\vect{C}$ with multiplicity 1. Then, there is a corresponding factor
\begin{equation}
q_j\left(s,\,P,\,D,\,\tau_1,\,\tau_2,\,\lambda_1=1\right)=s^2+Ds+P-\left(Dse^{-\tau_2s}+Pe^{-\tau_2s}\right)
\label{eq:centdir}
\end{equation}
in the characteristic quasi-polynomial \eqref{eq:chareq}. Without loss of generality, we will assign this eigenvalue to the state $\xi_1$. It can also be shown that the normalized eigenvector corresponding to $\xi_1$ is always $\vect{t}_1=1/\sqrt{n}\left[1\ 1\ \cdots\ 1\right]^T\in\R^n$, and it is selected as the first column of the earlier defined transformation matrix $\vect{T}$. We call this quantity, $\xi_1$, the weighted centroid, which is, obviously, topology dependent. It represents an agreement dynamics in the swarm. That is, if there is consensus we expect the behaviour of the weighted centroid to represent that agreement. The other factors of the characteristic equation in \eqref{eq:cefactdir} must then be related to some disagreement dynamics of the system. When these disagreement dynamics are stable (i.e., they vanish asymptotically) the agents should reach a consensus.

It is clear that $s=0$ is always a stationary root of \eqref{eq:centdir} independent of the delays, $\tau_1$ and $\tau_2$, which implies that the dynamics of the weighted centroid is at best marginally stable. In fact, provided that the disagreement dynamics is stable, $\xi_1$ reaches a non-zero steady state, when the topology is connected and has at least one spanning tree (a.k.a. arborescence) \cite{Gabow1978}. That is, each agent is connected to the root of its arborescence.  

If the communication topology does not have a spanning tree, however, $\lambda=1$ is a multiple eigenvalue of $\vect{C}$ and equation \eqref{eq:centdir} is a repeated factor within \eqref{eq:cefactdir}. These factors represent the dynamics of the centroids of the subgroups generated by the subgraphs which contain different spanning trees. If all the disagreement factors are stable the swarm members within a subgroup reach some stationary positions which are generally different. Thus the consensus is not achieved. These facts are stated in the following lemma. 

\begin{lemma}\emph{Group behaviour}.
\label{lemma}
Assume the communication topology has at least one spanning tree. Then, the agents in the group reach a consensus if and only if the factor \eqref{eq:centdir} is marginally stable and all the remaining factors of \eqref{eq:centdir} are stable. Furthermore, the consensus value will be $\sqrt{n}\xi_1\left(t\rightarrow\infty\right)$, whereas the other states $\xi_j\left(t\rightarrow\infty\right)=0$  for $j=2,\,3,\ldots,\,n$.
\end{lemma}

\emph{Proof:} First, we prove the necessity condition. From the definition of the state $\boldsymbol{\xi}$ we have $\left[\xi_1\left(t\right)\,\xi_2\left(t\right)\cdots\xi_n\left(t\right)\right]^T=\vect{T}\left[x_1\left(t\right)\,x_2\left(t\right)\cdots\,x_n\left(t\right)\right]$. If consensus is reached, the agents have a common steady state value which we denote by $\bar{x}=\lim_{t\rightarrow\infty}x_j\left(t\right)$, $j=1,\,2,\,3,\ldots,\,n$. Then:
\begin{equation}
\lim_{t\rightarrow\infty}\left[\xi_1\left(t\right)\,\xi_2\left(t\right)\,\cdots\,\xi_n\left(t\right)\right]^T=\bar{x}\vect{T}^{-1}\left[1\,1\,\cdots\,1\right]^T
\label{eq:limxi}
\end{equation}
Since the communication topology is assumed to have a spanning tree, $\lambda_1=1$ is a simple eigenvalue of $\vect{C}$, corresponding to the eigenvector $\vect{t}_1=1/\sqrt{n}\left[1\,1\,\cdots\,1\right]^T$, the first column of the transformation matrix $\vect{T}$. Since $\vect{T}^{-1}\left[1\,1\cdots\,1\right]^T=\sqrt{n}\vect{T}^{-1}\vect{t}_1=\sqrt{n}\left[1\,0\,\cdots\,0\right]^T$, equation \eqref{eq:limxi} leads to $\lim_{t\rightarrow\infty}\xi_1\left(t\right)=\sqrt{n}\bar{x}$, which indicates marginal stability for \eqref{eq:centdir} and $\lim_{t\rightarrow\infty}\xi_1\left(t\right)=0$ for $j=2,\,3,\ldots,\,n$ indicating asymptotic stability in the other factors of \eqref{eq:cefactdir}.

Next the sufficiency clause is proven. If \eqref{eq:centdir} is marginally stable (due to the characteristic root at the origin) and all the other factors in \eqref{eq:cefactdir} are stable, the steady state value of $\xi_1$ will be constant whereas the remaining $\xi_j\left(t\right)$, $j=2,\,3,\ldots,\,n$ tend to zero. Then,
$\lim_{t\rightarrow\infty}\left[\xi_1\left(t\right)\,\xi_2\left(t\right)\,\cdots\,\xi_n\left(t\right)\right]^T=\left[\sqrt{n}\bar{x}\,0\,\cdots\,0\right]^T$
Using the inverse transformation, the original states become:   
\begin{equation}
\lim_{t\rightarrow\infty}\left[x_1\left(t\right)\,x_2\left(t\right)\,\cdots\,x_n\left(t\right)\right]^T=\vect{T}\left[\bar{x}\,0\,\cdots\,0\right]^T=\bar{x}\vect{t}_1=\left[\bar{x}\,\bar{x}\,\cdots\,\bar{x}\right]^T
\label{eq:limx}
\end{equation}
implying the agents reach consensus.\hfill $\Box$

According to Lemma \ref{lemma}, a stable swarm consensus is reached if all the factors\newline$\factorm$, with $j=2,\,2,\,3,\ldots,\,\ell+m$, and the weighted centroid behaviour $q_1\left(s,\,P,\,D,\,\tau_1,\,\tau_2,\,\lambda=1\right)$ are all stable. In the contrary, if all the factors of \eqref{eq:centdir} are stable except $q_1$ in \eqref{eq:cefactdir}, i.e., the weighted centroid of the group has an unstable behaviour, the agents will move in coherence but along an unstable trajectory.
\section{CTCR Paradigm and the novel SDS domain}
\label{sec:III}
We just showed that the stability of the swarm dynamics is only achieved provided that each of the factors given in \eqref{eq:2ndfactor} and \eqref{eq:4thfactor} represent stable system behavior. For a given set of control parameters ($P$, $D$) these factors have the general formations of:
\begin{equation}
q\left(s,\,\tau_1,\,\tau_2\right)=g_{11}\left(s\right)+g_{12}\left(s\right)e^{-\tau_1\,s}+g_{13}\left(s\right)e^{-\tau_2\,s}
\label{eq:generic2ndfactor}
\end{equation}
for \eqref{eq:2ndfactor}, and
\begin{equation}
\begin{split}
q\left(s,\,\tau_1,\,\tau_2\right)=&g_{21}\left(s\right)+g_{22}\left(s\right)e^{-\tau_1\,s}+g_{23}\left(s\right)e^{-\tau_2\,s}+\\
&g_{24}\left(s\right)e^{-\left(\tau_1+\tau_2\right)\,s}+g_{25}\left(s\right)e^{-2\tau_2\,s}+g_{26}\left(s\right)e^{-2\tau_1\,s}
\end{split}
\label{eq:generic4thfactor}
\end{equation}
for \eqref{eq:4thfactor}.  In fact, by definitions of \eqref{eq:2ndfactor} and \eqref{eq:4thfactor}, the class of quasi-polynomials given in \eqref{eq:generic2ndfactor} is a subset of those in \eqref{eq:generic4thfactor}. Therefore we will concentrate on the stability treatment of the latter generic class. Notice that this quasi-polynomial represents a system with multiple and rationally independent time delays with delay cross-talk and commensurate degree 2.
 
The stability analysis of the factors in \eqref{eq:generic4thfactor} is not trivial, as the infinite dimensionality is doubled due to the presence of two rationally independent delays. This task is performed deploying a unique methodology called the Cluster Treatment of Characteristic Roots (CTCR) \cite{Olgac2005,Fazelinia2007,Ergenc2007}. The main philosophy behind it is the detection of the right half plane characteristic roots (i.e., unstable roots), of a Linear Time Invariant - Multiple Time Delay System (LTI-MTDS) (12). It is well known in the linear systems literature that the number of unstable roots can change only along certain loci (in this case in the domain of the delays). The CTCR method requires the exhaustive knowledge of these, so called, \emph{stability switching curves}. In this paper we utilize a novel approach to obtain this knowledge in a new domain which is called the Spectral Delay Space (SDS). The following paragraphs present some preparatory definitions and explain the key propositions of CTCR, in light of the works of \cite{Olgac2005} and \cite{Fazelinia2007}.
\begin{definition}
\label{def:kernel}
\emph{Kernel hypercurves} $\wp_0$: The hypercurves in the $\R^{2+}$ domain that consist exhaustively of all the points $\left(\tau_1,\,\tau_2\right)\in\R^{2+}$ which cause an imaginary root of \eqref{eq:generic4thfactor} at $s=\pm\omega i$ and satisfy the constraint $0<\tau_k\omega<2\pi$, $k=1,\,2$, are called the \emph{kernel hypercurves}. The points on these hypercurves contain the smallest possible delay values that create the given pair of imaginary roots at the frequency $\omega$.
\end{definition}
\begin{definition}
\label{def:offspring}
\emph{Offspring hypercurves} $\wp$: The hypercurves obtained from the kernel by the following pointwise non-linear transformation:
\begin{equation}
\left\langle\tau_1+\frac{2\pi}{\omega}j_1,\ \tau_2+\frac{2\pi}{\omega}j_2\right\rangle,\quad j_1,\,j_2=1,2,\ldots
\label{eq:nltran}
\end{equation}
are called the \emph{offspring hypercurves}. They are created by the periodicity of the imaginary roots with respect to the time delay.
\end{definition}

Jointly the kernel and the offspring hypercurves create all the loci in the delay space $\left(\tau_1, \tau_2\right)\in\mathbb{R}^2$, where the system in \eqref{eq:chareq} has imaginary characteristic roots. this complete set constitutes the departure point of the CTCR paradigm. It is crucial that an appropriate computational tool is used to capture all these hypercurves exhaustively. One should, however, not confuse this procedure with the CTCR paradigm itself, which uniquely declares the stability region exhaustively. for instance in a recent article, \cite{Sipahi2011} fall in this fallacy. They misrepresent a novel procedure of determining the kernel and offspring hyperplanes (they call the \emph{potential stability switching curves} instead) as a "new" stability declaration method, despite the fact that the CTCR is still the umbrella paradigm. This fallacy, indeed, results in the completely incorrect declarations of stability regions for the graphically displayed case study in \cite{Sipahi2011}, as well as misleading the readers on the novelties of the work.

We now present one more definition before the two key propositions of CTCR. 

\begin{definition}
\label{def:rt}
\emph{Root Tendency}, $RT$: At any point $\boldsymbol{\tau}\in\wp_0\cup\wp$ an infinitesimal increase in any of the individual delays, $\tau_j$, creates a transition of the root. Such transition can be to the right or to the left half of the complex plane. The Root Tendency, $RT$, indicates the direction of this transition as only one of the delays, $\tau_j$, increases by $\epsilon$, $0<\epsilon<<1$, while all the others remain constant:
\begin{equation}
\left.RT\right|_{s=\omega i}^{\tau_j}=\mathbf{sgn}\left[\Real\left(\left.\frac{\partial s}{\partial\tau_j}\right|_{s=\omega i}\right)\right]
\label{eq:rt}
\end{equation}
Clearly these root tendencies are $-1$ for stabilizing and $+1$ for destabilizing root crossings on the imaginary axis.
\end{definition}

\begin{proposition}
\label{p:prop1}
\emph{Finite Number of Kernel Hypercurves} A given LTI-MTDS can exhibit only a finite number, $m$, of kernel hypersurfaces. This number is upperbounded by the square of the order of the system: $m<n^2$. 
\end{proposition}
\begin{proposition}
\label{p:prop2}
\emph{Invariance of the Root Tendency} Take an imaginary characteristic root, $s=\omega i$, caused by any one of the infinitely many grid points on the kernel and offspring hypercurves in $\left(\tau_1,\,\tau_2\right)\in\R^{2+}$. The root tendency of these imaginary roots remains invariant so long as the grid points on different offspring hypersurfaces are obtained by keeping one of the delays fixed. That is, the root tendency with respect to the variations of $\tau_j$ is invariant from the kernel to the corresponding offspring, obtained by the non-linear mapping \eqref{eq:nltran}, as the other delay is fixed. 
\end{proposition}

\emph{Spectral Delay Space (SDS)}: We describe a new procedure in this segment for determining the kernel (and offspring) hypercurves. It is developed on a new domain called the Spectral Delay Space, SDS \cite{Fazelinia2007}. SDS is defined by the coordinates $\nu_j=\tau_j\omega$ for every point $\left(\tau_1,\,\tau_2\right)\in\R^{2+}$ on the kernel or the offspring hypercurves. This is a conditional mapping: if a delay set $\left(\tau_1,\,\tau_2\right)$ creates an imaginary root $\omega i$, (i.e., if the point is on the kernel or the offspring hypercurves) then  $\left(\tau_1\omega,\,\tau_2\omega\right)$ forms a point in the SDS.  On the contrary, $\left(\tau_1,\,\tau_2\right)$ points that do not generate an imaginary root have no representation in the SDS domain.

The main advantage of SDS is that the representation of the kernel hypercurve in the SDS, denoted as $\wp_0^{SDS}$, which is named the \emph{building hypercurve}, is confined into a square of edge length $2\pi$ (see Definition \ref{def:kernel}). Then, it is only necessary to explore a finite domain in SDS in order to find the representation of the building hypercurves in SDS. This finite domain is known as the building block (BB), i.e., a square of $2\pi\times2\pi$ for the case of two rationally independent delays. Similarly we name the corresponding representation of offspring hypercurves in the SDS as the \emph{reflection curves}. An important advantage of SDS domain is that the transitions from the building hypercurves to the reflection hypercurves is achieved simply by translating the building hypercurve by $\pi$2, as opposed to using the pointwise non-linear transformation \eqref{eq:nltran}. We prevent the undesirable shape distortion from kernel to offspring which occurs in the delay space. In the contrary, the reflection curves in the SDS, $\wp^{sds}$, are generated by simply stacking the BB squares (with $wp_0^{SDS}$ in it) one on top of each other. Naturally, this property is named the \emph{stackability feature} of the SDS. There are several other intriguing properties of the SDS and BB concepts which can be found in the works of \cite{Fazelinia2007} and \cite{FazeliniaThesis}.

With these definitions and propositions, we perform the preparatory step of CTCR. It is the exhaustive determination of all the stability switching curves, for the generic factors $\factorm$ as in \eqref{eq:factorsdir} within the semi-infinite quadrant of  $\left(\tau_1,\,\tau_2\right)\in\R^{2+}$.

There are different ways in which the characteristic equation can be transformed to the SDS domain to obtain the building hypercurves. The following mathematical procedure, which uses a half-angle tangent substitution, is taken from the appendix of the thesis by \cite{FazeliniaThesis}.

We start by replacing $s=\omega i$ into \eqref{eq:generic4thfactor}, and then the exponential terms are replaced by:
\begin{equation}
e^{\tau_k\omega\,i}=\cos\left(\nu_k\right)+i\,\sin\left(\nu_k\right),\quad \nu_k=\tau_k\omega,\quad k=1,2
\label{eq:exp}
\end{equation}
Then the sine and cosine functions are expressed in terms of a half-angle tangent substitution:
\begin{equation}
\cos\left(\nu_k\right)=\frac{1-z_k^2}{1+z_k^2},\quad \sin\left(\nu_k\right)=\frac{2z_k}{1+z_k^2},\quad z_k=\tan\left(\nu_k\right)
\label{eq:halftan}
\end{equation}
Now \eqref{eq:generic4thfactor} can be written as a polynomial in $\omega$ with complex coefficients $c_k$ parametrized in $z_1$ and $z_2$:
\begin{equation}
q_j\left(\omega,z_1,z_2\right)=\sum_{k=0}^2{c_k\left(P,\,D,\,\lambda_j,\,z_1,\,z_2\right)\left(\omega\,i\right)^k}=0
\label{eq:factorz}
\end{equation}
If there is a solution $\omega\in\R^+$ to \eqref{eq:factorz}, both its real and imaginary parts must be zero simultaneously:
\begin{subequations}
\label{eq:reandim}
\begin{align}
\Re\mathfrak{e}\left(q_j\left(\omega,z_1,z_2\right)\right)&=\sum_{k=0}^2{f_k\left(z_1,\,z_2\right)\omega^k}=0\label{eq:re}\\
\Im\mathfrak{m}\left(q_j\left(\omega,z_1,z_2\right)\right)&=\sum_{k=0}^2{g_k\left(z_1,\,z_2\right)\omega^k}=0\label{eq:im}
\end{align}
\end{subequations}
The condition for \eqref{eq:reandim} to have a common root is simply stated using a Sylvester's matrix:
\begin{equation}
\vect{M}=\left[\begin{array}{cccc}
f_2\left(z_1,\,z_2\right)&f_1\left(z_1,\,z_2\right)&f_0\left(z_1,\,z_2\right)&0\\
0&f_2\left(z_1,\,z_2\right)&f_1\left(z_1,\,z_2\right)&f_0\left(z_1,\,z_2\right)\\
g_2\left(z_1,\,z_2\right)&g_1\left(z_1,\,z_2\right)&g_0\left(z_1,\,z_2\right)&0\\
0&g_2\left(z_1,\,z_2\right)&g_1\left(z_1,\,z_2\right)&g_0\left(z_1,\,z_2\right)
\end{array}\right]
\end{equation}
In order to simultaneously satisfy bot \eqref{eq:re} and \eqref{eq:im}, $\vect{M}$ should be singular. This results in the following expression in terms of $z_1$ and $z_2$, or $\nu_1$ and $\nu_2$:
\begin{equation}
\det\left(\vect{M}\right)=F\left(z_1,\,z_2\right)=F\left(\tan\left(\nu_1\right),\,\tan\left(\nu_2\right)\right)
\label{eq:detsylv}
\end{equation}
which constitutes a closed form description of the kernel curves in the SDS $\left(\nu_1,\,\nu_2\right)$, i.e., the Building Block curves. To obtain its graphical depiction, one of the parameters, say $\nu_2$, can be scanned in the range of $\left[0,\,2\pi\right]$ and the corresponding $\nu_1$ values are calculated again in the same range. Notice that every point $\left(\nu_1,\,\nu_2\right)$ on these curves brings an imaginary characteristic root at $\pm\omega\,i$. That is, we have a continuous sequence of $\left(\nu_1,\,\nu_2,\,\omega\right)$ sets all along the building block curves. 

In order to assess the stability properties of the system in the space of the delay, we now back transform the building and reflection hypercurves from  the $\left(\nu_1,\,\nu_2\right)$ domain of the SDS to the $\left(\tau_1,\,\tau_2\right)$ delay space, using the inverse transformation of \eqref{eq:exp} and generate the kernel and offspring hypercurves. The application of the root tendency invariance property follows, resulting in the complete and exact stability outlook of the system, which is unique. Examples of this construction are presented in the following section.

\section{Illustrative Examples}
\begin{figure}[!tb]
\center
\includegraphics[scale=0.8]{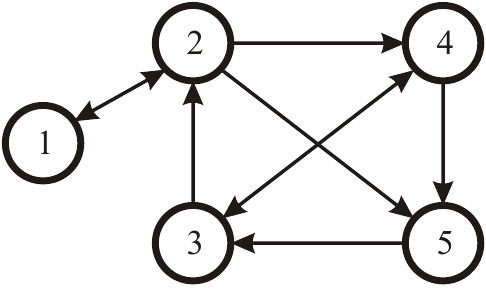}
\caption{Example directed communication topology with 5 agents.}
\label{fig:top5dir}
\end{figure}

In this section we present the simulation results that verify the theoretical analysis of previous sections. All of the examples use the communication topology of Fig.~\ref{fig:top5dir}. The corresponding weighted adjacency matrix is:
\begin{equation}
\vect{C}=\left[\begin{array}{ccccc}
0&1&0&0&0\\
0.5&0&0&0.5&0\\
0&0.5&0&0.5&0\\
0&0&0.5&0&0.5\\
0&0.5&0.5&0&0
\end{array}\right]
\end{equation}
and its eigenvalues are are $1$, $0.38$, $-0.5$, and $-0.44\pm0.37i$. The control gains used in the examples are $P=2$ and $D=0.8$.

\begin{figure}[!tb]
\center
\includegraphics[scale=0.3]{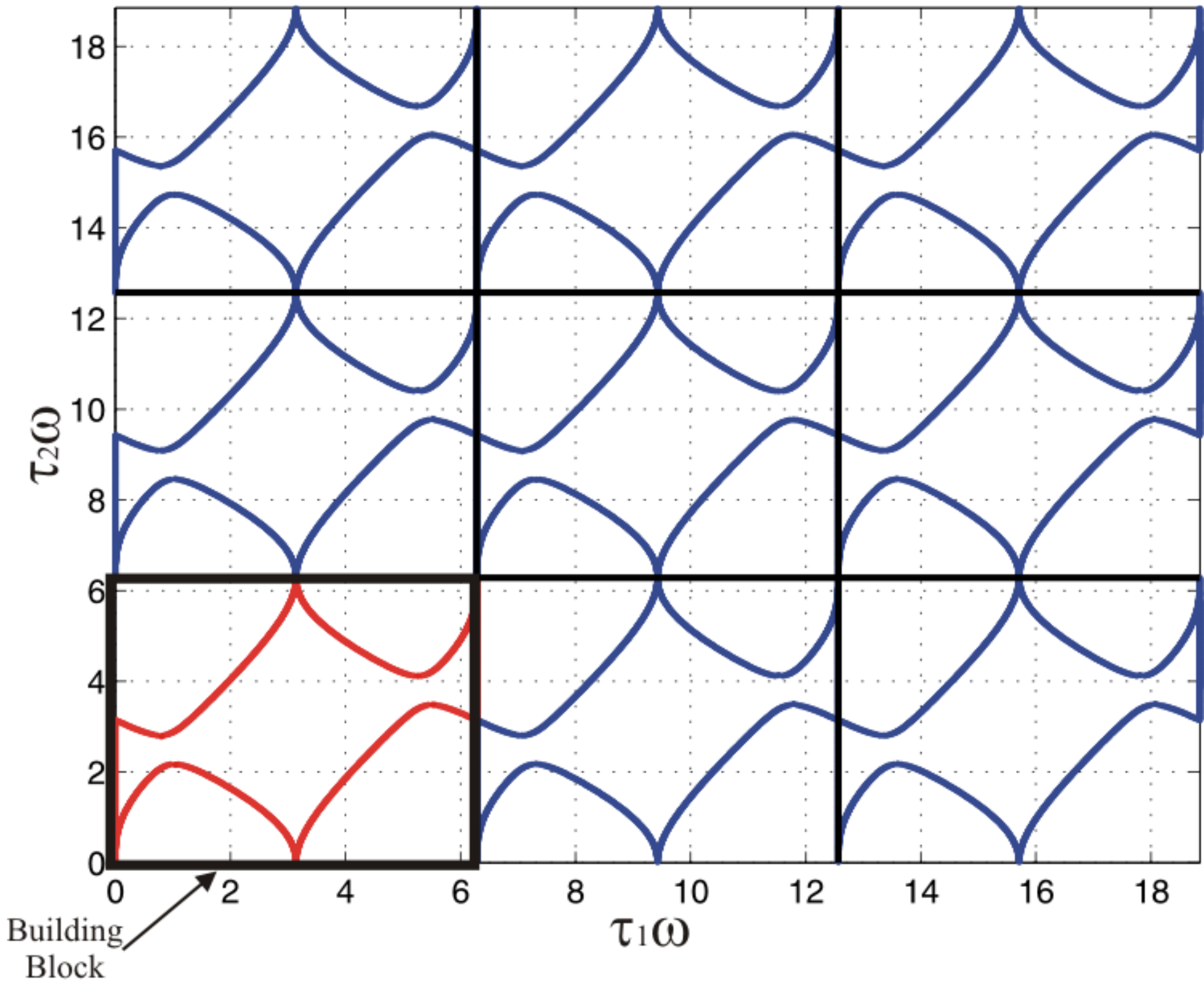}
\caption{Delay Space representation of the stability switching curves generated by the factor corresponding to $\lambda=1$ with $P=2$ and $D=0.8$.}
\label{fig:bb1}
\end{figure}

\begin{figure}[!tb]
\center
\includegraphics[scale=0.7]{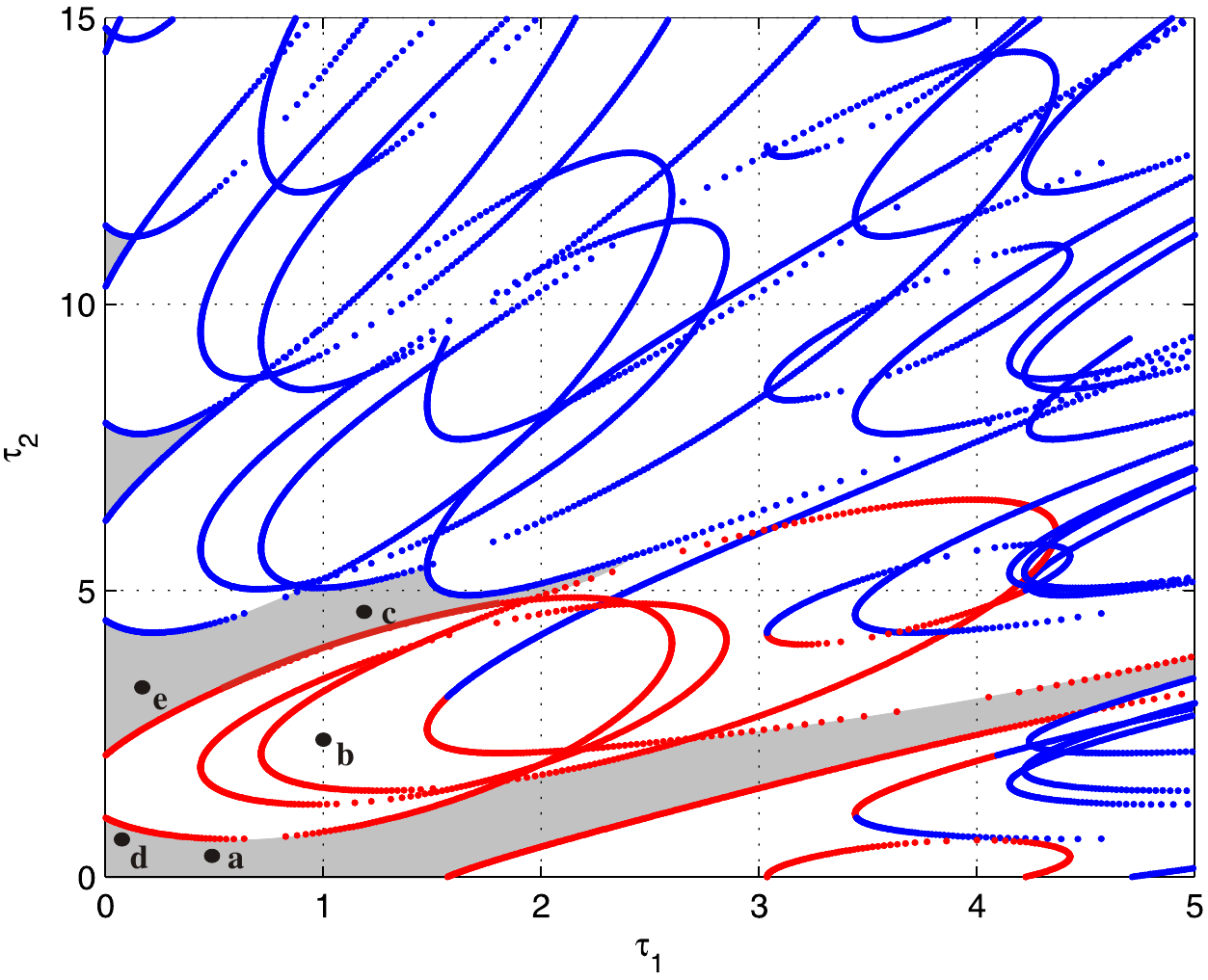}
\caption{Complete stability picture for the communication topology of Fig.~\ref{fig:top5dir}, using $P=2$ and $D=0.8$.}
\label{fig:regdir}
\end{figure}

The deployment of the CTCR paradigm using the building block concept leads to Figs.~\ref{fig:bb1} and \ref{fig:regdir}. Figure \ref{fig:bb1} shows the spectral delay space (SDS) representation of the characteristic factor corresponding to $\lambda=1$, with the respective building and reflection hypercurves. Other factors are excluded to prevent congestion. The properties of the SDS are clear: the reflection curves, presented in blue, are obtained by translating the building curves, which are depicted in red. Notice that the building curves are confined to a square of side length $2\pi$ i.e., the building block.
 
\begin{figure}[!tb]
\center
\includegraphics[scale=0.4]{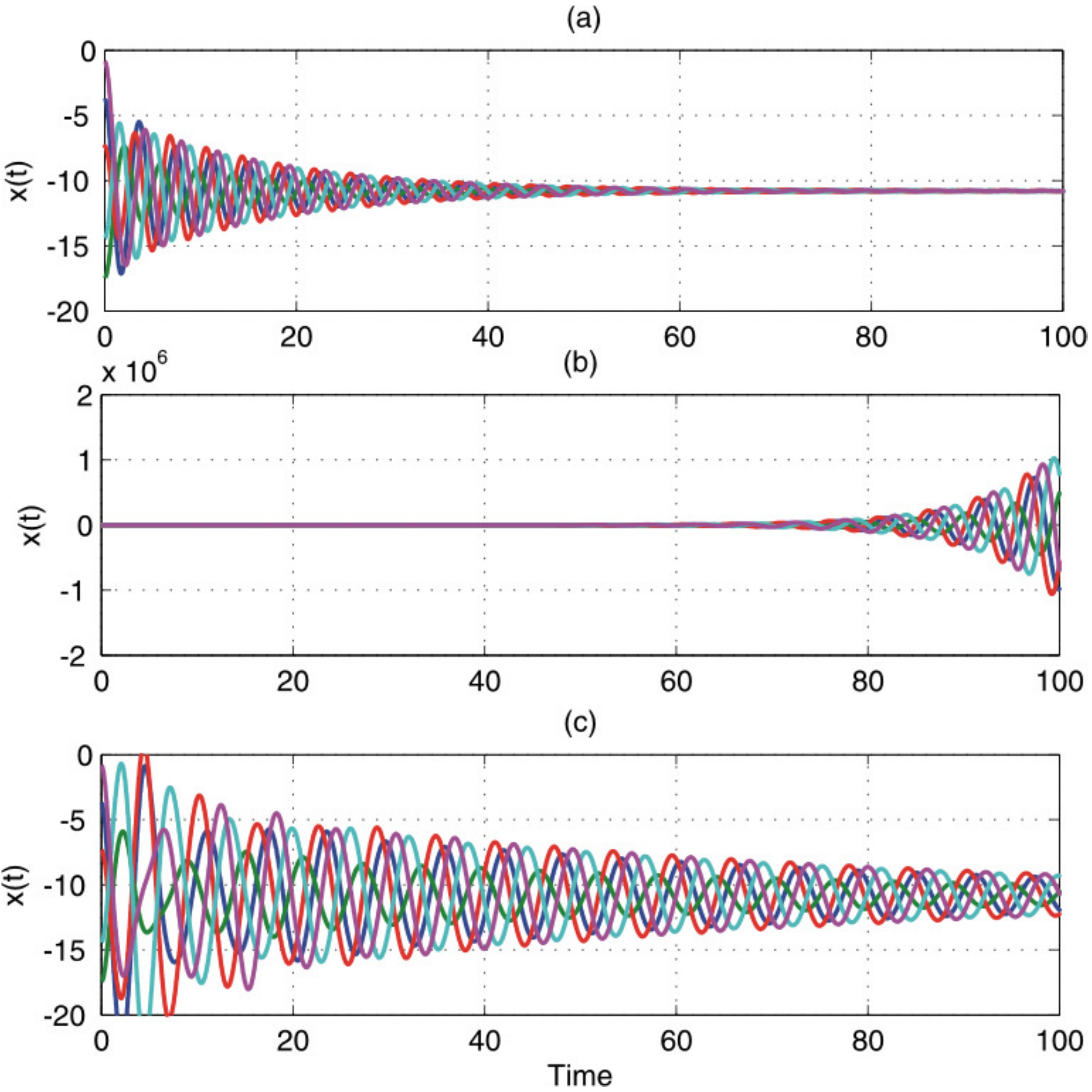}
\caption{Example agents' behavior for three different points in Fig.~\ref{fig:regdir}. From top to bottom, points \textbf{a}, \textbf{b} and \textbf{c}.}
\label{fig:simulation}
\end{figure}

The superposition of the stability regions of the individual factors results in the stability chart of Fig.~\ref{fig:regdir}. Delay combinations inside the shaded portion of this figure result in stable behavior, i.e., agents reach consensus. These results are verified in Fig.~\ref{fig:simulation}, in which panel (a) shows the simulation results for a delay combination of $\left(\tau_1,\,\tau_2\right)=\left(0.5,\,0.5\right)$ seconds, corresponding to point \textbf{a} in Fig.~\ref{fig:regdir}, panel (b) is for $\left(\tau_1,\,\tau_2\right)=\left(1,\,2.5\right)$ (point \textbf{b} in Fig.~\ref{fig:regdir}) and panel (c) is for $\left(\tau_1,\,\tau_2\right)=\left(1.3,\,4.5\right)$, represented by point \textbf{c} in Fig.~\ref{fig:regdir}. It is clear that the delay combinations corresponding to the points inside the shaded region (\textbf{a} and \textbf{c}) generate stable consensus, whereas points outside create divergent behavior.

\subsection{Delay Scheduling Concept}
A more interesting, and counter-intuitive, conclusion we can draw from Figs.~\ref{fig:regdir} and \ref{fig:simulation} is that larger time delays not always imply worsening control performance of the system. Both time delays are larger in \textbf{c} than in \textbf{b}, however, the system is stable in the first case and unstable in the second. This feature has been used in a control synthesis scheme named the \emph{delay scheduling} \cite{Olgac2011}. It offers wonderful options for the most effective control parameter selections, because we can manipulate the feedback delay as a part of the control design. Due to the causality, the delay can only be increased (i.e., one cannot sense the events before they become available) and the control designer can, however, intentionally prolong a given set of delays further. This is a paradoxical observation. Below we present a brief discussion on how one should perform the selection of increase in delay, a.k.a. the delay scheduling, along with a deployment scheme of this concept concept in group dynamics.

We first establish the stability chart of the particular topology following the steps described in the earlier sections (i.e., Fig.~\ref{fig:regdir}). Then the present inherently delay composition is monitored in real time. Let us assume these delays render the operating point of \textbf{b} in Fig.~\ref{fig:regdir} which is in the unstable region. The controller can now impose an intentional extra delay, as a remedy, to bring the combined delay values to point \textbf{c}, where the stability is recovered. It is obvious that such a point must be selected within the stable region. How to select the point \textbf{c} in the stable operating region is, however, still an open question. In order to answer this question one needs a crisp stability picture, such as Fig.~\ref{fig:regdir}, for which CTCR is the only procedure available today.

The challenge of determining where to take the operating point by introducing extra delays (i.e., the choice of point \textbf{c}) offers an additional freedom to the control designer. We aim, here, to select a stable point which facilitates the fastest disturbance rejection capabilities to the swarm formation.  As such we pose the operational question of the speed for reaching consensus. This information is crucial in tuning the delays (i.e., prolonging) in order to achieve the most desirable (expeditious) speed in reaching consensus.  We explain the relevant steps next on how the stability map, which is given in Fig.~\ref{fig:regdir}, can be utilized for this purpose.

\begin{figure}[!tb]
\center
\includegraphics[scale=0.15]{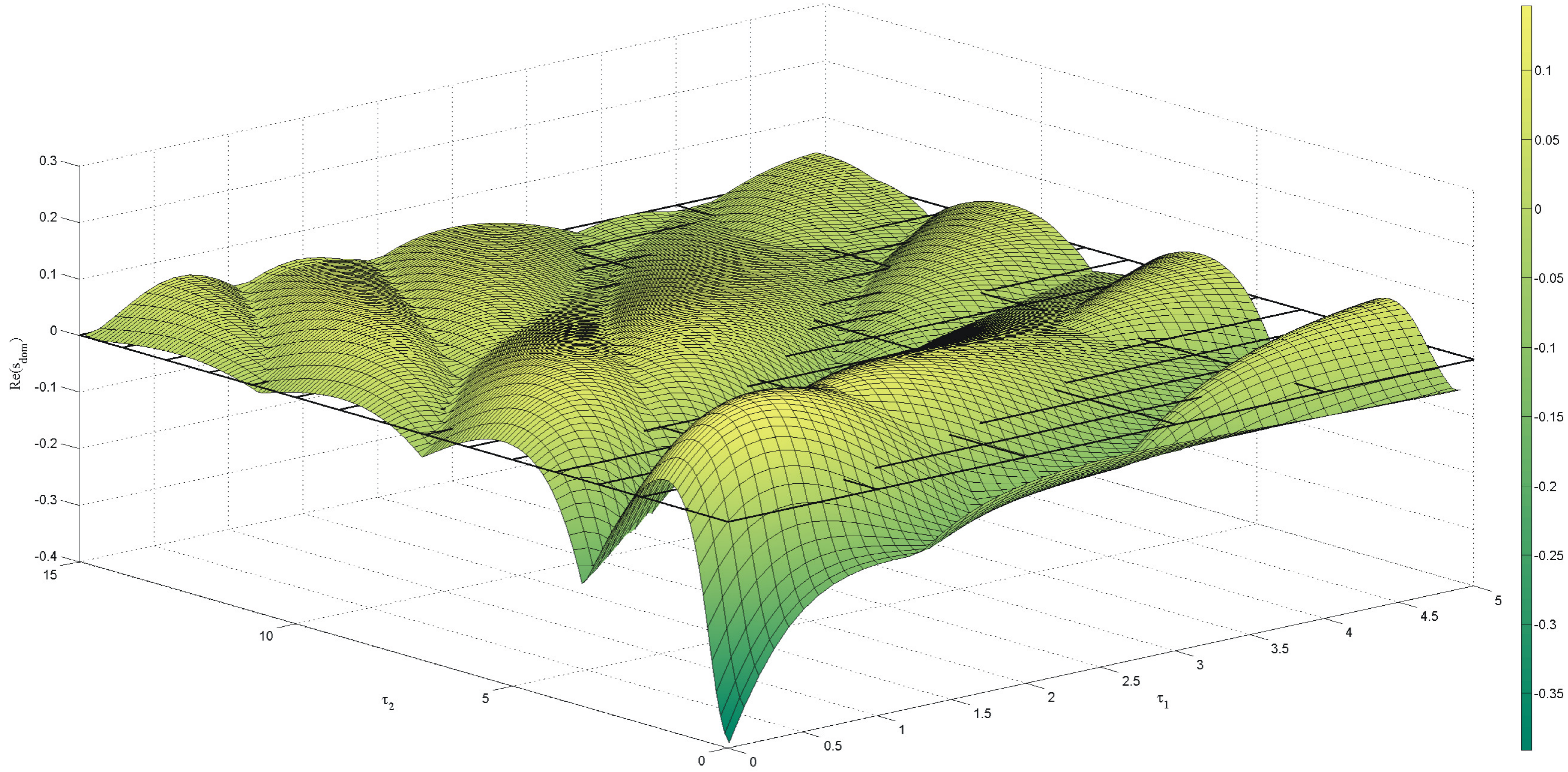}
\caption{Variation of the real part of the dominant root for the complete system. Grid represents the zero real part partitioning.}
\label{fig:domroot}
\end{figure}

The speed of group consensus is dictated by the dominant characteristic root of the infinite dimensional system given in \eqref{eq:c2delay2}.  The dominant root of the characteristic equation \eqref{eq:chareq} is the rightmost of the dominant roots of all the factors given. These factors, $\factorm$, are quasi-polynomials and their roots can only be determined by some numerical approximations, like those proposed by \cite{vyhlidal2009} and \cite{Breda2006}. We deploy the QPmR (Quasi-Polynomial, mapping-based Root finding) algorithm of \cite{vyhlidal2009} here. For each point $\left(\tau_1,\,\tau_2\right)$ we determine the dominant root and its real part, $\Real\left(s_{dom}\right)$. For the specific structure in this example case study, it is displayed in Fig.~\ref{fig:domroot}.

We wish to stress several observations over this figure:
\begin{enumerate}
\item The stable operating zones of Fig.~\ref{fig:regdir} are approximately generated by the point-wise evaluation of the dominant roots over a grid in $\left(\tau_1,\,\tau_2\right)$ space. Obviously they are the regions where $\Real\left(s_{dom}\right)<0$.
\item The time constants expected in reaching consensus, $-1/\Real\left(s_{dom}\right)$, in example case \textbf{a} of Fig.~\ref{fig:regdir} concurs with the location of $s_{dom}$. For this point $\Real\left(s_{dom}\right)=-0.0610$, which corresponds to a time constant of 16.3934 seconds, and a settling time (the speed of reaching consensus) of approximately 66s, as observed in Fig.~\ref{fig:simulation}.
\item Although the fastest consensus is achieved at $\tau_1=\tau_2=0$, larger delays do not always mean longer consensus time. Compare for instance points \textbf{d} and \textbf{e} in Fig.~\ref{fig:regdir}. For the first one, $\left(\tau_1,\,\tau_2\right)=\left(0.05,\,0.8\right)$ and $\Real\left(s_{dom}\right)=-0.04$, whereas for \textbf{e} $\left(\tau_1,\,\tau_2\right)=\left(0.1,\,3.5\right)$ and $\Real\left(s_{dom}\right)=-0.05$. Obviously, the latter makes a better choice from the point of speed in reaching consensus, although both delays are larger in the second case. This is the core concept of delay scheduling [26-28]. 
\end{enumerate}

In order to trigger the delay selections under this methodology, one has to have the complete system stability tableau in the delay space, that is, Figs.\ref{fig:regdir} and \ref{fig:domroot}. The procedures which are discussed in this paper to obtain the stability outlook, offer a very efficient technique. For instance, to create one of the individual stability tableaus for a given factor in (8) we encounter an average computational cost of 1s of CPU time, including display operations, on a 2.93 GHz Intel Core 2 Duo - based computer with 4 GB RAM, running Matlab 2010a without an elaborate code optimization for speed.
\section{Conclusions} In this paper we study the consensus problem for multi-agent systems, considering directed communication topologies and multiple time delays. Two constant and rationally independent time delays are taken into account, one affecting the position and the other the velocity information exchange among the agents.  Both delays can be large and they are assumed to be rationally independent.

The procedure starts with an enabling factorization of the characteristic equation of the system. The factors involved are only distinguished from one another by a set of eigenvalues of a particular system matrix. The process simplifies the stability analysis considerably and converts it into a repetitive application of a reduced-order core stability analysis n times. The Cluster Treatment of Characteristic Roots (CTCR) methodology and the new concept of Spectral Delay Space are utilized to obtain the exact stability bounds in the domain of the time delays. As the number of agents increase, the only procedural complexity is in determining the set of eigenvalues of a given system matrix. Numerical example cases show the validity of the results.

We also look at the variations on the speed of reaching consensus as the two delays in the communication channels vary. It becomes clear that one may have the option to intentionally prolong the delays to achieve a faster consensus (a procedure called the Delay Scheduling). This ability becomes available only if precise knowledge of the stability picture in the parametric space is available. That is what the CTCR paradigm uniquely provides.

\bibliographystyle{plain}
\bibliography{IJC3}
\end{document}